\documentclass[twoside]{article}
\usepackage{amsmath,amsthm,amssymb,stmaryrd,latexsym,amsbsy,xy}
\usepackage{amscd, amssymb, verbatim}
\usepackage{pstricks,pst-plot,pst-node}
\usepackage{amsfonts}
\usepackage{anysize}
\usepackage[pdftex,bookmarks=true]{hyperref}
\usepackage{fancyhdr}

\marginsize{2.9cm}{2.9cm}{2.9cm}{2.9cm}

\newtheorem{theorem}{Theorem}[section]

\newtheorem{lemma}[theorem]{Lemma}
\newtheorem{prop}[theorem]{Proposition}
\newtheorem{cor/conj}[theorem]{Corollary/Conjecture}
\newtheorem{cor}[theorem]{Corollary}

\newtheorem{conj/cor}[theorem]{Conjecture/Corollary}

\numberwithin{equation}{section}

\newcommand{\bpf}{\noindent {\em Proof.  }}
\newcommand{\epf}{\qed \vspace{+10pt}}

\newcommand{\Sym}{\mathrm{Sym}}
\newcommand{\AT}{A_{\bT}^*}

\newcommand{\ATorb}{A_{\bT, \mathrm{orb}}^*}

\newcommand{\vir}{\mathrm{vir}}

\newcommand{\Spec}{\mathrm{Spec}}

\newcommand{\M}{\overline{M}}
\newcommand{\Mc}{\overline{M}^{\circ}}
\newcommand{\I}{\overline{I}}
\newcommand{\age}{\mathrm{age}}
\newcommand{\Aut}{\mathrm{Aut}}
\newcommand{\conn}{\mathrm{conn}}

\newcommand{\ev}{\mathrm{ev}}
\newcommand{\spa}{\textrm{ }}
\newcommand{\Aone}{\mathcal{A}_1}
\newcommand{\Ar}{\mathcal{A}_{r}}
\newcommand{\An}{\mathcal{A}_{n}}

\newcommand{\tC}{\tilde{C}}

\newcommand{\tf}{\tilde{f}}
\newcommand{\tSi}{\tilde{\Si}}
\newcommand{\Eij}{\cE_{ij}}
\newcommand{\ttt}{(t_1+t_2)}

\newcommand{\fB}{\mathfrak{B}}
\newcommand{\fS}{\mathfrak{S}}

\newcommand{\cB}{\mathcal{B}}
\newcommand{\cC}{\mathcal{C}}
\newcommand{\cE}{\mathcal{E}}

\newcommand{\cL}{\mathcal{L}}
\newcommand{\cP}{\mathcal{P}}
\newcommand{\cQ}{\mathcal{Q}}

\newcommand{\bC}{\mathbb{C}}

\newcommand{\bP}{\mathbb{P}}
\newcommand{\bQ}{\mathbb{Q}}
\newcommand{\bT}{\mathbb{T}}
\newcommand{\bZ}{\mathbb{Z}}

\newcommand{\al}{\alpha}
\newcommand{\be}{\beta}
\newcommand{\g}{\gamma}

\newcommand{\de}{\delta}
\newcommand{\si}{\sigma}
\newcommand{\Si}{\Sigma}
\newcommand{\la}{\lambda}

\newcommand{\Zc}{\mathrm{Z}^\circ}
\newcommand{\wt}{\widetilde}

\makeatletter
\renewcommand*\@fnsymbol[1]{\the#1}

\begin{document}

\title{\Large \bf Connected Gromov-Witten invariants of $[\Sym^n(\Ar)]$}
\author{
Wan Keng Cheong\footnote{Department of Mathematics, California Institute of Technology, Pasadena, CA 91125, USA. \newline Email address: keng@caltech.edu} \spa and 
Amin Gholampour\footnote{Department of Mathematics, California Institute of Technology, Pasadena, CA 91125, USA. \newline Email address: agholamp@caltech.edu}
}

\date{}
\maketitle

\setcounter{section}{-1}

\begin{abstract}
We explore the theory of connected Gromov-Witten invariants of the symmetric product stack $[\Sym^n(\Ar)]$. We derive closed-form expressions for all equivariant invariants with two insertions and reveal a natural correspondence between the theory and the relative Gromov-Witten theory of $\Ar\times \bP^1$. When $n\leq 3$, we determine 3-point (usual) Gromov-Witten invariants of $[\Sym^n(\Aone)]$.
\end{abstract}

\section{Overview} 
Let $\Ar$ be the minimal resolution of the type $A_r$ surface singularity and $[\Sym^n(\Ar)]$ the $n$-th symmetric product stack of $\Ar$. These spaces admit $\bT$-actions, where $\bT=(\bC^\times)^2$ is a torus.

In this article, we study the $\bT$-equivariant, genus $0$, connected Gromov-Witten theory of the stack $[\Sym^n(\Ar)]$, which is defined over the space of twisted stable maps where the covering associated to the domain is connected. The connected invariants with two insertions are completely solved here. In fact, we have the following.

\begin{theorem}\label{main0}
Two-point extended connected invariants of $[\Sym^n(\Ar)]$ of nonzero degree can be expressed by one-part double Hurwitz numbers, and so they admit closed-form expressions. 
\end{theorem}

This result recovers certain relative Gromov-Witten invariants of $\Ar \times \bP^1$ (c.f. Proposition \ref{relative stable}) and suggests a natural correspondence between the connected Gromov-Witten theory of symmetric product stacks and the relative Gromov-Witten theory of threefolds. We remark that it also plays a special role in the study of orbifold Gromov-Witten theory of $[\Sym^n(\Ar)]$ (c.f. \cite{C2}).

This paper is organized as follows. In Section \ref{s:1}, we recall some basic facts about the surface $\Ar$ and the orbifold cohomology for a symmetric product.
In Section \ref{s:2}, after introducing the theory of connected Gromov-Witten invariants of $[\Sym^n(\Ar)]$, we establish Theorem \ref{main0}.
We also calculate $3$-point usual Gromov-Witten invariants of $[\Sym^n(\Aone)]$ for $n \leq 3$.

Additionally, the first author would like to take this opportunity to thank his adviser Tom Graber for many stimulating discussions and guidance.

\paragraph{Convention} We need cohomology and homology groups with rational or integral coefficients. To avoid doubling indices, we identify $A^i(X)=H^{2i}(X; \bQ)$, $A_i(X; \bZ)=H_{2i}(X; \bZ)$ and so on, for any complex variety $X$ to appear in this paper. We call them cohomology or homology groups rather than Chow groups.

\section{Preliminaries}\label{s:1}
\paragraph{$A_r$-resolutions}
The resolved surface $\Ar$ is a smooth toric surface. Under the $\bT$-action, it has $r+1$ fixed points $x_1,\ldots, x_{r+1}$, ordered so that $x_i$ and $x_{i+1}$ are connected by a $\bT$-fixed curve $E_i$. Indeed, the chain of curves $\cup_{i=1}^r E_i$ is exactly the exceptional locus of $\Ar$. The intersection matrix $(E_i \cdot E_j)$ is minus the Cartan matrix for the $A_r$ Dynkin diagram.
Note that the set $\{E_1, \ldots, E_r\}$ is a basis for $A_1(\Ar; \bZ)$. Its dual basis $\{\omega_1,\ldots, \omega_r\}$ will also be considered later on.

The equivariant cohomology $A_{\bT}^*(\mathrm{point})$ of a point is $\bQ[t_1,t_2],$ with $(t_1, t_2)$ the tangent weights at $0\in \bC^2$.
We let
$$
(L_i, R_i)=((r-i+2)t_1+(1-i)t_2, \spa (-r+i-1)t_1+it_2),
$$
which are actually the tangent weights at $x_i$.

\paragraph{Orbifold cohomology and a basis}
Given any smooth complex variety $X$, we consider the $n$-th symmetric product $\Sym^n(X):=X^n/\fS_n$, where the symmetric group $\fS_n$ acts on $X^n$ by permuting coordinates. The space $\Sym^n(X)$ is generally singular but admits a smooth stack (orbifold) structure given by the quotient stack $[\Sym^n(X)]$. Suppose the torus $\bT$ acts on $X$, by $\bT$-equivariant orbifold cohomology $$\ATorb([\Sym^n(X)])$$ of the orbifold $[\Sym^n(X)]$, we mean the $\bT$-equivariant cohomology $\AT(\I[\Sym^n(X)])$ of the stack $\I[\Sym^n(X)]$ of cyclotomic gerbes in $[\Sym^n(X)]$. Also, the connected components of $\I[\Sym^n(X)]$ are in one-to-one correspondence with the partitions of $n$. For each partition $\la$ of $n$, the corresponding component is $[\overline{X(\la)}]:=X^g/\overline{C(g)}$, for $g$ a permutation of cycle type $\la$, $X^g$ the $g$-fixed locus of $X$ and $\overline{C(g)}$ the quotient group of the centralizer of $g$ by $\langle g \rangle$.

Let $X=\Ar$ for the rest of this section. As seen in Section 2 of \cite{C2}, the cohomology-weighted partition $\la(\vec{\eta})$, partition $\la$ with parts labelled with classes $\eta_i$'s on $\Ar$, defines a class on the component $[\overline{X(\la)}]$. We denote the class by $\la(\vec{\eta})$ as well. Note that the classes $\la(\vec{\eta})$'s, running over all partitions $\la$ of $n$ and entries $\eta_i \in \{ 1, E_1, \ldots, E_r \}$, give a basis for $\ATorb([\Sym^n(X)])$. By replacing $1$, $E_1, \ldots, E_r$ with $\bT$-fixed point classes $[x_1], \ldots, [x_{r+1}]$, we obtain a basis, referred to as fixed-point basis, for the localized cohomology $\ATorb([\Sym^n(X)])\otimes_{\bQ[t_1, t_2]}\bQ(t_1, t_2)$.

Let's define the age of $\la(\vec{\eta})$, denote by $\age(\la(\vec{\eta}))$, to be the age of the component $[\overline{X(\la)}]$, i.e. $\age(\la):=n-\ell(\la)$. Also, we need some more notation for later convenience. Let
$$(2)=1(1)^{n-2}2(1), \spa \spa D_k= 1(1)^{n-1}1(\omega_k), \spa k=1,\ldots, r,$$
which form a basis for divisors on $[\Sym^n(X)]$.
We also use the symbol 
$(2)$
to stand for the partition $(1^{n-2}, 2)$ of length $n-1$.

\section{Connected Gromov-Witten theory} \label{s:2}
\subsection{Definitions}\label{s: 2-1}
\paragraph{Twisted stable maps} 
We assume the reader is familiar with the notions of orbifold Gromov-Witten theory \cite{CR2, AGV}.
In what follows, we will restrict ourselves to $X=\Ar$ or $X=\Spec \spa \bC$ and identify $A_1(\Sym^n(X); \bZ)$ with $A_1(X; \bZ).$ Also, we let $\cB G$ be the classifying stack of $G$, i.e. $[\Spec \spa \bC/G]$, for any finite group $G$.

Given any effective curve class $\be \in A_1(X; \bZ)$, the moduli space 
$$\overline{M}_{0,k}([\Sym^n(X)], \beta)$$ 
parametrizes $k$-pointed genus $0$ twisted stable map $f:(\cC,\cP_1,\ldots, \cP_k)\to [\Sym^n(X)],$
where markings $\cP_1,\ldots, \cP_k$ and nodes of $\cC$ are twisted points. In particular, $\cP_i\cong \cB \mu_{r_i}$ for some cyclic group $\mu_{r_i}$ of order $r_i$. The map $f$ comes with an ordinary stable map $f_c:(C, c(\cP_1),\ldots, c(\cP_k)) \to \Sym^n(X)$ of degree $\be$. Here $C$ is the coarse moduli space of $\cC$ and $c$ is the canonical map to the coarse moduli space.
 
We denote by
$
\ev_i: \overline{M}_{0,k}([\Sym^n(X)], \beta) \to \I[\Sym^n(X)]
$
the evaluation morphism at the $i$-th marking. For partitions $\si_1,...,\si_k$ of $n$, let 
$$\overline{M}([\Sym^n(X)],\si_1,...,\si_k; (a,\beta))=\bigcap_{i=1}^k \ev_i^{-1}([\overline{X(\si_i)}])\cap \bigcap_{j=1}^a \ev_{k+j}^{-1}([\overline{X((2))}])$$ 
be a substack of the moduli space $\overline{M}_{0,k+a}([\Sym^n(X)], \beta)$. The space possesses a virtual fundamental class $[\overline{M}([\Sym^n(X)],\si_1,...,\si_k; (a,\beta))]^{\vir}$ of dimension $n \cdot \dim(X)+k-3-\sum_{i=1}^k \age(\si_i)$.
For simplicity, put
$$
\overline{M}(\cB\fS_n,\si_1,...,\si_k; a):=\overline{M}(\cB\fS_n,\si_1,...,\si_k; (a,0)).
$$

Let $f:\cC, (\cP_1, \ldots, \cP_k; \cQ_1, \ldots, \cQ_a) \to  [\Sym^n(X)]$ be a twisted stable map representing an element of the moduli space $\M([\Sym^n(X)], \si_1, \ldots, \si_k;(a,\be))$, the markings $\cP_1, \ldots, \cP_k$ are referred to as distinguished marked points and the other markings $\cQ_1, \ldots, \cQ_a$ are referred to as simple marked points. The map $f$ gives rise to a Cartesian diagram
\begin{equation}
\begin{CD} 
\cC^\prime @> f^\prime >> X^n\\
@V VV @VV V\\
\cC @>f >> [\Sym^n(X)] \\
\end{CD}\label{twisted}
\end{equation}
which induces a morphism $\tf:\tC \to X$ by taking $f^\prime$ mod $\fS_{n-1}$ and composing with the $n$-th projection. The composition $\tC \to \cC \overset{c}{\rightarrow} C$ is a degree $n$ admissible cover (see \cite{ACV}), branched with monodromy $\si_1, \ldots, \si_k, (2), \ldots, (2)$ above the $k+a$ marked points.

Now let
$$\Mc([\Sym^n(X)], \si_1, \ldots, \si_k;(a,\be))$$ 
be the component of $\M([\Sym^n(X)], \si_1, \ldots, \si_k;(a,\be))$ which parametrizes connected covers (i.e., each induced cover $\tC$ is connected). 

\paragraph{$k$-point invariants}
Given cohomology classes $\al_i \in \ATorb([\Sym^n(\Ar)])$ ($i=1,\ldots, k$), the $k$-point, extended, $\bT$-equivariant Gromov-Witten invariant of $[\Sym^n(\Ar)]$ is defined by
$$
\left\langle  \al_1, \ldots, \al_k  \right\rangle_{(a,\be)}^{[\Sym^n(\Ar)]}=\frac{1}{a!} \sum_{|\si_1|, \ldots, |\si_k|=n} \int_{[\M([\Sym^n(\Ar)], \si_1, \ldots, \si_k;(a,\be))]^{\vir}_{\bT}}\ev_1^*(\al_1) \cdots \ev_k^*(\al_k),
$$
where the symbol $[\spa \spa ]^\vir_{\bT}$ stands for the $\bT$-equivariant virtual class. Note that the (non-extended) Gromov-Witten invariant $\left\langle  \al_1, \ldots, \al_k  \right\rangle_{\be}^{[\Sym^n(\Ar)]}$ is simply $\left\langle  \al_1, \ldots, \al_k  \right\rangle_{(0,\be)}^{[\Sym^n(\Ar)]}$.

We will not calculate these invariants but will mention some related facts later. In this article, we would like to focus on the $k$-point, extended, connected invariant defined by
$$
\left\langle  \al_1, \ldots, \al_k  \right\rangle_{(a,\be)}^{[\Sym^n(\Ar)], \conn}=\frac{1}{a!} \sum_{|\si_1|, \ldots, |\si_k|=n}  \int_{[\Mc([\Sym^n(\Ar)], \si_1, \ldots, \si_k;(a,\be))]^{\vir}_{\bT}}\ev_1^*(\al_1) \cdots \ev_k^*(\al_k).
$$
To distinguish the connected Gromov-Witten invariant from the usual one, we add a superscript `conn'. The space $\Mc([\Sym^n(\Ar)], \si_1, \ldots, \si_k;(a,\be))$ is compact whenever $\be \ne 0$; hence the corresponding invariant takes value in $\bQ[t_1,t_2]$.

We will drop the superscript $[\Sym^n(\Ar)]$ if there is no ambiguity. Also, we say that an extended invariant is of nonzero (resp. zero) degree if the associated 2-tuple $(a,\be)$ has $\be \ne 0$ (resp. $\be=0$). 

\subsection{Closed formulas for 2-point extended invariants}
\paragraph{Double Hurwitz numbers}
For partitions $\la, \rho$ of $n$, the double Hurwitz number $H^g_{\la,\rho}$ is the (weighted) number of possibly disconnected genus $g$ branched covers of $\bP^1$ with ramification profiles $\la, \rho, \eta_1, \ldots, \eta_b$ for some $b$. Here all $\eta_i$'s are $(2)$ and $b=2g-2+\ell(\la)+\ell(\rho)$ is determined by the Riemann-Hurwitz formula. Note that we don't label ramification points.

In general, it is not easy to obtain a closed formula for $H^g_{\la,\rho}$. However, when $\rho=(n)$, we have the following fact due to Goulden, Jackson, and Vakil.
\begin{prop}[\cite{GJV}] \label{GJV}
Given any partition $\la=(\la_1, \ldots, \la_{\ell(\la)})$ of $n$.
The double Hurwitz number $H_{\la,(n)}^{g}$ is the coefficient of $t^{2g}$ in the power series expansion of
$$
\frac{(2g+\ell(\la)-1)! \spa n^{2g+\ell(\la)-2}}{|\Aut(\la)|}{\frac{t/2}{sinh (t/2)}} \prod_{i=1}^{\ell(\la)}\frac{sinh (\la_i t/2)}{\la_i t/2}.
$$
\end{prop}

\paragraph{Localization calculation}
We aim to calculate all 2-point, extended, connected invariants of nonzero degree. To achieve this, we only have to study
\begin{equation} \label{extended invariant}
\left\langle \mu_1(\vec{\g}_1),\mu_2(\vec{\g}_2) \right\rangle_{(a,\be)}^{\conn},
\end{equation}  
for partitions $\mu_1$, $\mu_2$ of $n$, nonnegative integer $a$, nonzero effective curve class $\be \in A_1(X; \bZ)$, and cohomology classes $\g_{1k}$, $\g_{2\ell}$'s which are $E_1, \ldots, E_r$ or $1$. 
For the rest of the section, we fix such $n$, $a$, $\mu_1$, and $\mu_2$.

Given nonnegative integers $a_1$, $a_2$, let $g_{a_1}$, $g_{a_2}$, $g(a)$ be integers satisfying 
$a_k=2g_{a_k}-1+\ell(\mu_k)$ for $k=1,2$, and $g(a)=\frac{1}{2}(a-\ell(\mu_1)-\ell(\mu_2)+2)$. Put $\Eij=E_i+\cdots+E_j$ for $i\leq j$ and $[\mu_k(\vec{\g}_k)]=|\Aut(\mu_k)|/|\Aut(\mu_k(\vec{\g}_k))|$ for $k=1,2$.

Our extended invariants can be expressed in terms of double Hurwitz numbers.
\begin{theorem}\label{main}
Assume that $\g_{1k}$, $\g_{2\ell}$'s are $E_1, \ldots, E_r$ or $1$ and $\be$ is a nonzero effective curve class. If $\be=d\Eij$ for some $d,i,j$ and all $\g_{1k}$, $\g_{2\ell}$'s are either $E_i$ or $E_j$, the invariant \eqref{extended invariant} is given by
\begin{equation}\label{Hurwitz}
\frac{\ttt(-1)^{g(a)}(-1-\de_{1,r})^{\ell(\mu_1)+\ell(\mu_2)} d^{a-1}[\mu_1(\vec{\g}_1)] [\mu_2(\vec{\g}_2)]}{n^{a-2}} \sum_{a_1+a_2=a} \frac{H_{\mu_1,(n)}^{g_{a_1}}H_{\mu_2,(n)}^{g_{a_2}}}{a_1!a_2!},
\end{equation}
where $\de_{1,r}$ is the Kronecker delta (which is $1$ if $r=1$ and $0$ otherwise).
Otherwise, the invariant \eqref{extended invariant} vanishes.
\end{theorem}
\bpf
Let $r>1$. Each insertion of \eqref{extended invariant} is a linear combination of fixed point classes with coefficients zero or having nonnegative $\ttt$-valuation. According to Proposition 4.1 in \cite{C2}, \eqref{extended invariant} is divisible by $t_1+t_2$. 

As mentioned earlier, \eqref{extended invariant} is a polynomial in $t_1$, $t_2$. If at least one of $\g_{1k}$, $\g_{2\ell}$'s is 1, \eqref{extended invariant} must be zero because the sum of the degrees of the insertions is at most $\ell(\mu_1)+\ell(\mu_2)-1$, which is the virtual dimension.

Assume that all $\g_{1k}$, $\g_{2\ell}$'s are $E_1, \ldots, E_r$. In light of the dimension and $\ttt$-divisibility, \eqref{extended invariant} is proportional to $t_1+t_2$. By the proof of Proposition 4.1 in \cite{C2}, \eqref{extended invariant} is zero if $\be$ is not a multiple of $\Eij$ for all $i,j$. Now we consider $\be=d\Eij$ for some $d,i,j$. As we may evaluate \eqref{extended invariant} modulo $\ttt^2$, any $\bT$-fixed locus that contributes a factor $\ttt^k$ for some $k\geq 2$ may be ruled out. That is, it is enough to investigate those $\bT$-fixed loci $\M_{a_1,a_2}$'s with the following configuration:

Let $[f]$ be an element in $\M_{a_1,a_2}$, $a_1+a_2=a$, with induced map $\tf$. The twisted curve $\cC$ decomposes into three pieces, namely $\cC_{a_1}\cup \Si \cup \cC_{a_2}$: for $k=1,2$, $\cC_{a_k}$ is a contracted component carrying $a_k$ simple markings, and its associated cover $\tC_{a_k}$ is of genus $g_{a_k}$; the intersection $\cC_{a_1}\cap \cC_{a_2}$ is empty; and $\Si$ is a chain of non-contracted components, which connects $\cC_{a_1}$ and $\cC_{a_2}$. Note that $\cC_{a_k}$'s are twisted points whenever they contain less than three special points and are otherwise twisted curves. 

However, in order for the contribution of $\M_{a_1,a_2}$ to \eqref{extended invariant} not to vanish, the ramification points lying above the distinguished markings must map to $x_i$ or $x_{j+1}$. As a result \eqref{extended invariant} vanishes if one of $\g_{1k}$, $\g_{2\ell}$'s is $E_k$ for some $k \ne i,j$. This completes the proof of the second assertion.

Now we show the first assertion. Assume that all $\g_{1k}$, $\g_{2\ell}$'s are $E_i$ or $E_j$. 
By the second assertion, the invariant \eqref{extended invariant} is congruent modulo $\ttt^2$ to
\begin{equation}\label{replace} 
[\mu_1(\vec{\g}_1)] [\mu_2(\vec{\g}_2)] \left\langle \mu_1(-\frac{1}{L_i}[x_i], \ldots, -\frac{1}{L_i}[x_i]),\mu_2(-\frac{1}{R_{j+1}}[x_{j+1}], \ldots, -\frac{1}{R_{j+1}}[x_{j+1}]) \right\rangle_{(a,d\Eij)}^{\conn}. 
\end{equation}

We may use \eqref{replace} to make a further reduction: only those $\bT$-fixed loci, denoted by $F_{a_1,a_2}$'s, with an additional property that the marking corresponding to $\mu_1$ is in $\cC_{a_1}$ and the one corresponding to $\mu_2$ is in $\cC_{a_2}$ can make contributions. Keep in mind that the component $F_{a_1,a_2}$ must contribute a single factor of $t_1+t_2$, in which case the associated covering $\tSi$ of $\Si$ is a chain of rational curves totally branched over two points (which must be either nodes or markings) and branched nowhere else.
Therefore we reduce our calculation to the integral over
$$
\overline{M}(\cB\fS_n, \mu_1,(n); \spa a_1)\times \overline{M}(\cB\fS_n, \mu_2, (n) ; a_2)
$$
followed by division by the product of the automorphism factor $d^{j-i+1}$ and the distribution factor $a_1! a_2!$ of simple marked points. 

Let $\epsilon_1:\overline{M}(\cB\fS_n, \mu_1,(n); \spa a_1) \to \overline{M}_{0,a_1+2}$ be the natural morphism mapping $\cC_{a_1}$ to its coarse moduli space $C_{a_1}$ (the node $\cC_{a_1}\cap \Si$ is mapped to the marking $Q_1$) and $\cL_1$ the tautological line bundle formed by the cotangent space $T^*_{Q_1} C_{a_1}$. Let $\psi_1=c_1(\cL_1)$.
We define $\epsilon_2:\overline{M}(\cB\fS_n, \mu_2,(n); \spa a_2) \to \overline{M}_{0,a_1+2}$ and $\psi_2$ in a similar way.
To proceed, let's summarize the contributions of virtual normal bundles. Set $\theta=(r+1)t_1$.
\begin{itemize}
\item Contracted components: For $k=1,2$, $\cC_{a_k}$ contributes
$$
(-1)^{g_{a_k}-1}\theta^{2g_{a_k}-2} \mod \ttt.
$$

\item A chain of non-contracted components: The contribution of each node-smoothing is just 
$(\frac{t_1+t_2}{d})^{-1}$. All other node contributions are
$L_k R_k, \spa \spa k=i,\cdots, j+1$, each of which equals $-\theta^2$ mod $t_1+t_2$. Furthermore, all non-contracted curves contribute
$(\frac{t_1+t_2}{-\theta^2})^{j-i+1} \mod \ttt^2$. Hence the total contribution equals
$$-\theta^2 d^{j-i}\ttt \mod \ttt^2.$$

\item Smoothing nodes joining a contracted curve to a non-contracted curve:
The contributions are given by
$$
\frac{1}{\frac{1}{n}(\frac{n R_i}{d}- \epsilon_1^* \psi_1)}, \spa \frac{1}{\frac{1}{n}(\frac{n L_{j+1}}{d}- \epsilon_2^* \psi_2)}.
$$
\end{itemize}
The contribution $I_{a_1,a_2}$ of the fixed locus $F_{a_1,a_2}$ to \eqref{replace} is congruent modulo $\ttt^2$ to
\begin{eqnarray*}
& & -\theta^2 d^{j-i}\ttt \spa \frac{[\mu_1(\vec{\g}_1)] [\mu_2(\vec{\g}_2)]}{d^{j-i+1} a_1!a_2!}\spa \theta^{\ell(\mu_1)}(-\theta)^{\ell(\mu_2)} \cdot \frac{(-1)^{a_1}}{\theta^4}\\ 
&\times& \prod_{k=1}^2 (-1)^{g_{a_k}} \theta^{2g_{a_k}} \frac{1}{n^{a_k-1}} \spa (\frac{d}{\theta})^{a_k} \spa \int_{\overline{M}(\cB\fS_n, \mu_k,(n); \spa a_k)}\epsilon_k^* \psi_k^{a_k-1}.
\end{eqnarray*}
In the first line of the above expression, e.g. the class $\theta^{\ell(\mu_1)}$ is congruent modulo $\ttt$ to the restriction of $\ev_1^{*}(\mu_1(-\frac{1}{L_i}[x_i], \ldots, -\frac{1}{L_i}[x_i]))$ to the component $F_{a_1,a_2}$. In addition, each factor in the second line is replaced with $1$ in case $a_k=0$.
Simplifying the expression yields
$$
\frac{\ttt(-1)^{g(a)+\ell(\mu_1)+\ell(\mu_2)}d^{a-1}[\mu_1(\vec{\g})] [\mu_2(\vec{\g}_2)]}{n^{a-2}a_1!a_2!}\int_{\overline{M}(\cB\fS_n, \mu_1,(n); a_1)}\epsilon_1^* \psi_1^{a_1-1} \spa \int_{\overline{M}(\cB\fS_n, \mu_2, (n) ; a_2)}\epsilon_2^* \psi_2^{a_2-1}.
$$
For $a_k>0$,
$$
\int_{\overline{M}(\cB\fS_n, \mu_k,(n); \spa a_k)}\epsilon_k^* \psi_k^{a_k-1}= \deg(\epsilon_k)\int_{\overline{M}_{0,a_k+2}} \psi_k^{a_k-1}=H_{\mu_k,(n)}^{g_{a_k}}.$$
We conclude that 
$$
I_{a_1,a_2}\equiv \frac{\ttt(-1)^{g(a)+\ell(\mu_1)+\ell(\mu_2)}d^{a-1}[\mu_1(\vec{\g}_1)] [\mu_2(\vec{\g}_2)]}{n^{a-2}} \cdot \frac{H_{\mu_1,(n)}^{g_{a_1}} H_{\mu_2,(n)}^{g_{a_2}}}{a_1!a_2!} \mod \ttt^2.
$$
Finally by \eqref{replace}, and keeping in mind that \eqref{extended invariant} is a multiple of $t_1+t_2$, we obtain \eqref{Hurwitz}.

The case where $r=1$ can be argued similarly and so is omitted.
\epf

Theorem \ref{main}, together with Proposition \ref{GJV}, shows Theorem \ref{main0}. By applying the intersection matrix with respect to the curve classes $E_1,\ldots, E_r$, we arrive at the following statement which is a little shorter than Theorem \ref{main}.
\begin{cor}\label{two point}
Let $\g_{1k}$, $\g_{2\ell}$'s be 1 or divisors on $\Ar$ and $\be$ a nonzero effective curve class. If $\be=d\Eij$ for some $d,i,j$, the connected invariant 
$\left\langle \mu_1(\vec{\g}_1),\mu_2(\vec{\g}_2) \right\rangle_{(a,\be)}^{\conn}$
is given by
$$
\frac{\ttt(-1)^{g(a)} d^{a-1}[\mu_1(\vec{\g}_1)] [\mu_2(\vec{\g}_2)]\prod_{k=1}^{\ell(\mu_1)}(\Eij \cdot \g_{1k}) \prod_{k=1}^{\ell(\mu_2)}(\Eij\cdot \de_k)}{n^{a-2}} \sum_{a_1+a_2=a} \frac{H_{\mu_1,(n)}^{g_{a_1}}H_{\mu_2,(n)}^{g_{a_2}}}{a_1!a_2!}.
$$
Otherwise, it is zero.
\end{cor}

\paragraph{Comparison to the relative theory}
Connected Gromov-Witten invariants of $[\Sym^n(\Ar)]$ can be related to certain relative invariants of the Gromov-Witten theory of $\Ar\times \bP^1$. To make this precise, we make some definitions.

The 3-point function of the connected Gromov-Witten theory of $[\Sym^n(\Ar)]$ encodes 3-point extended invariants:
$$\langle \langle \al_1, \al_2, \al_3 \rangle \rangle^{\conn}:=\sum_{a=0}^\infty \sum_{\be \in A_1(\Ar;\bZ)} \left\langle \al_1, \al_2, \al_3 \right\rangle_{(a,\be)}^{\conn}u^a s_1^{\be \cdot \omega_1}\cdots s_r^{\be \cdot \omega_r}.$$
On the other hand, we denote by
$$
Z^\circ (\Ar \times \bP^1)_{\la(\vec{\eta_1}), \rho(\vec{\eta_2}), \si(\vec{\eta_3})}\in \bQ(t_1,t_2)((u))[[s_1, \ldots, s_r]]
$$
the partition function of relative Gromov-Witten theory of the threefold $\Ar \times \bP^1$ with connected domain curves of arbitrary genus and with relative conditions given by cohomology-weighted partitions $\la(\vec{\eta_1}), \rho(\vec{\eta_2}), \si(\vec{\eta_3})$. For more information, see \cite{M}, \cite{C2}.

The degree $0$ connected invariants are in connection with the relative connected invariants:
\begin{lemma}\label{deg zero}
Let $\mu_1(\vec{\g}_1)$ and $\mu_2(\vec{\g}_2)$ be as in Corollary \ref{two point}. For $\theta=1(1)^n$ or $(2)$, 
$$
\sum_{a=1}^\infty \left\langle \mu_1(\vec{\g}_1), \theta, \mu_2(\vec{\g}_2)\right\rangle_{(a,0)}^{\conn, [\Sym^n(\Ar)]}u^a=u^{\ell(\mu_1)+\ell(\mu_2)-\age(\theta)} \Zc (\Ar \times \bP^1)_{\mu_1(\vec{\g}_1), \theta, \mu_2(\vec{\g}_2)}|_{s_1=\cdots=s_r=0}.
$$
\end{lemma}
\bpf
Let $U_i\cong \bC^2$ be an open subset of $\Ar$ with $E_{i-1}$ and $E_i$ being $x$-axis and $y$-axis respectively. Since the associated cover is connected and collapses to fixed points of $\Ar$, we may write the left side of the above equation as 
$$\sum_{i=1}^s a_i \sum_{a=1}^\infty \left\langle \wt{\si_1}, \theta, \wt{\si_2} \right\rangle_{(a,0)}^{\conn, [\Sym^n(U_i)]}u^a,$$ 
where $\wt{\si_k}$'s are partitions labelled with fixed point classes on $\Ar$ and $a_i \in \bQ(t_1,t_2)$. Similarly, the right side can be expressed as 
$$u^{\ell(\mu_1)+\ell(\mu_2)-\age(\theta)}\sum_{i=1}^s b_i \spa \Zc (U_i \times \bP^1)_{\wt{\si_1}, \theta, \wt{\si_2}}|_{s_1=\cdots=s_r=0}.$$ 
It is easy to check that $a_i=b_i$ for each $i$, and therefore the lemma follows from the corresponding statement in the affine plane case (c.f. \cite{BP}, \cite{BG}).
\epf

The proof of the lemma also shows that two-point extended connected invariants of zero degree are determined by the Gromov-Witten theory of $[\Sym^n(\bC^2)]$. In other words, we have provided a complete solution to the 2-point connected theory. Further, we have the following correspondence claimed earlier.
\begin{prop}\label{relative stable}
Let $\mu_1(\vec{\g}_1)$ and $\mu_2(\vec{\g}_2)$ be as in Corollary \ref{two point}. For $\theta=1(1)^n$, $(2)$ or $D_k$, $k=1, \ldots, r$,
$$
\langle \langle \mu_1(\vec{\g}_1), \theta, \mu_2(\vec{\g}_2)\rangle \rangle^{\conn}=u^{\ell(\mu_1)+\ell(\mu_2)-\age(\theta)} \Zc (\Ar \times \bP^1)_{\mu_1(\vec{\g}_1), \theta, \mu_2(\vec{\g}_2)},
$$
\end{prop}
\bpf
From our determination of 3-point functions $\langle \langle \mu_1(\vec{\g}_1), \al, \mu_2(\vec{\g}_2)\rangle \rangle^{\conn}$ and Proposition 4.3 in \cite{M}, the proposition follows from Lemma \ref{deg zero}. 
\epf

\paragraph{The case $n \leq 3$ and $r=1$}
For any $\al_1$, $\al_2$, $\al_3\in \AT([\Sym^n(\Aone)])$, we define the 3-point function of $[\Sym^n(\Aone)]$ by Gromov-Witten invariants:
$$
\langle \langle  \al_1, \al_2, \al_3 \rangle \rangle ^{[\Sym^n(\Aone)]}=\sum_{d=0}^\infty \langle \al_1, \al_2, \al_3 \rangle^{[\Sym^n(\Aone)]}_{(0, dE_1)}s_1^d.
$$
Clearly, the case of $n=1$ is determined by Theorem \ref{main} and the classical orbifold cup product. In general, we need the following fact about 2-point extended Gromov-Witten invariants (not necessarily connected) of nonzero degree.
\begin{prop}[\cite{C2}] \label{2-point extended}
Given classes $\la(\vec{\eta_1})$, $\rho(\vec{\eta_2})$ of $\ATorb([\Sym^n(\Ar)])$ where the entries of $\vec{\eta_1}$ and $\vec{\eta_2}$ are $1$ or divisors on $\Ar$. For any curve class $\be \ne 0$, the identity
$$
\left\langle \la(\vec{\eta_1}), \rho(\vec{\eta_2}) \right\rangle_{(a,\be)}=\sum \langle \si({\vec{\xi_1}}) , \si({\vec{\xi_2}})\rangle \left\langle \mu_1(\vec{\g_1}), \mu_2(\vec{\g_2})\right\rangle_{(a,\be)}^{\conn}
$$
holds. Here $\langle \si({\vec{\xi_1}}) , \si({\vec{\xi_2}})\rangle$'s are orbifold pairings and the sum is taken over all possible classes $\si({\vec{\xi}_1})$, $\si({\vec{\xi}_2})$, $\mu_1(\vec{\g_1})$ and $\mu_2(\vec{\g_2})$ such that $\la(\vec{\eta_1})=\si({\vec{\xi}_1})\mu_1(\vec{\g_1})$ and $\rho(\vec{\eta_2})=\si({\vec{\xi}_2})\mu_2(\vec{\g_2})$.
\epf
\end{prop}

This actually determines 2-point extended invariants of nonzero degree (c.f. \cite{C2}). As an application, we can deduce the following consequence for $3$-point functions of $[\Sym^n(\Aone)]$.
\begin{prop} \label{CRC 3}
Let $n \leq 3$. For any classes $\al_1, \al_2, \al_3 \in \ATorb([\Sym^n(\Aone)])$, the 3-point function 
\begin{equation}\label{specialization}
\langle \al_1, \al_2, \al_3 \rangle^{[\Sym^n(\Aone)]}
\end{equation}
is completely determined.
\end{prop}
\bpf
We only exhibit the proof in case $n=3$, the case of $n=2$ being similar and much easier. 
Let's consider the following orthogonal basis $\fB$ for $\ATorb([\Sym^3(\Aone)])$. Observe that the elements are grouped according to their ages.
\begin{eqnarray*}
&a_0=1(1)^3, \spa a_1=1(1)^2 1(E_1), \spa a_2=1(1)1(E_1)^2, \spa a_3=1(E_1)^3& \\
&b_0=1(1) 2(1), \spa b_1= 1(E_1) 2(1) , \spa b_2=1(1)2(E_1), \spa b_3=1(E_1)2(E_1)&\\
&c_0=3(1), \spa c_1=3(E_1).&
\end{eqnarray*} 

By linearity, we only have to examine 3-point function \eqref{specialization} for $\al_1, \al_2, \al_3 \in \fB$. 
First of all, we have two immediate facts.
\begin{enumerate}
\item[(1)] 3-point functions with one insertion from $a_0$, $a_1$ or $b_0$ are known.

\item[(2)] Up to the ordering of the insertions, 3-point functions other than $\langle a_i, a_j, a_k \rangle$, $\langle b_i, b_j, a_k \rangle$, $\langle c_i, c_j, a_k \rangle$, $\langle b_i, b_j, c_k \rangle$, $\langle c_i, c_j, c_k \rangle$ are identically zero for monodromy reasons.
\end{enumerate}

We now focus on those possibly nonzero 3-point functions. Their determination can be divided into the following steps. (For ease of explanation, we suppress the superscript $[\Sym^3(\Aone)]$ and create an extra symbol: for any 3-point functions $f,g$, we say that $f\sim g$ if $f-g$ has been determined).
\begin{enumerate}
\item[(3)] $\langle a_i, a_j, a_k \rangle $: For $k=2$, by the WDVV equation and (1), (2),
$$
\langle a_i, a_j, a_2 \rangle \langle a_2^\vee, a_1, a_1 \rangle + \langle a_i, a_j, a_3 \rangle \langle a_3^\vee, a_1, a_1 \rangle =\sum_{\al \in \fB} \langle a_i, a_1, \al \rangle \langle \al^\vee, a_j, a_1 \rangle \sim 0.
$$
One check straightforwardly that $\langle a_3^\vee, a_1, a_1 \rangle_0=0$. By Proposition \ref{2-point extended} and Theorem \ref{main}, $\langle a_3^\vee, a_1, a_1 \rangle_{dE_1}=0$ for each $d>0$. Further, $\langle a_2^\vee, a_1, a_1 \rangle_0$ is nonzero and so the inverse of $\langle a_2^\vee, a_1, a_1 \rangle$ exists. Thus $\langle a_i, a_j, a_2 \rangle \sim 0$. 
Now we are left with the case $\langle a_3, a_3, a_3 \rangle$.
By WDVV and (1), (2) again, $\langle a_3, a_3, a_3 \rangle \langle a_3^\vee, a_1, a_2 \rangle$ is determined by
$$
\langle a_3, a_1, a_2 \rangle \langle a_2^\vee, a_3, a_2 \rangle + \langle a_3, a_1, a_3 \rangle \langle a_3^\vee, a_3, a_2 \rangle - \langle a_3, a_3, a_2 \rangle \langle a_2^\vee, a_1, a_2 \rangle \sim 0.
$$
We deduce that $\langle a_3, a_3, a_3 \rangle \sim 0$ by non-vanishing of $\langle a_3^\vee, a_1, a_2 \rangle_0$.

\item[(4)] The proof of $\langle b_i, b_j, a_k \rangle \sim 0$ and $\langle c_i, c_j, a_k \rangle \sim 0$ proceeds exactly as in (3).

\item[(5)] $\langle \tau_1, \tau_2, c_0 \rangle $ for $\tau_1,\tau_2=b_i$ or $c_j$: By (1)-(4) and WDVV, 
$$
\langle \tau_1, \tau_2, c_0 \rangle \langle c_0^\vee, b_0, b_0 \rangle + \langle \tau_1, \tau_2, c_1 \rangle \langle c_1^\vee, b_0, b_0 \rangle \sim 0.
$$
Thus $\langle \tau_1, \tau_2, c_0 \rangle \sim 0$ since $\langle c_1^\vee, b_0, b_0 \rangle=0$ (by a straightforward check on the degree $0$ term and by Theorem \ref{main} on all other terms) and $\langle c_0^\vee, b_0, b_0 \rangle_0 \ne 0$.
\item[(6)] $\langle \tau_1, \tau_2, c_1 \rangle $ for $\tau_1,\tau_2=b_i$ or $c_1$: 
$$
\langle \tau_1, \tau_2, c_1 \rangle \sim \langle c_1^\vee, c_0, a_1 \rangle^{-1} (\sum_{\al \in \fB} \langle \tau_1, c_0, \al \rangle \langle \al^\vee, \tau_2, a_1 \rangle - \langle \tau_1, \tau_2, c_0 \rangle \langle c_0^\vee, c_0, a_1 \rangle) \sim 0.
$$
\end{enumerate}

Thus our assertion follows from (1)-(6). 
\epf

\phantomsection

\end{document}